\documentclass[reqno]{amsart}
\usepackage{amssymb}
\usepackage{graphicx}

\usepackage[usenames, dvipsnames]{color}
\usepackage{verbatim}
\usepackage{mathrsfs}
\usepackage{bm}
\usepackage{cite}

\numberwithin{equation}{section}

\newtheorem{theorem}{Theorem}[section]

\theoremstyle{definition}
\newtheorem{remark}[theorem]{Remark}

\theoremstyle{definition}

\theoremstyle{definition}

\makeatletter
\def\dashint{\operatorname%
{\,\,\text{\bf-}\kern-.98em\DOTSI\intop\ilimits@\!\!}}
\makeatother

\def\\det{\text{\det}}

\def\.5{\frac{1}{2}}

\newcommand{\RN}[1]{%
  \textup{\uppercase\expandafter{\romannumeral#1}}%
}

\renewcommand{\epsilon}{\varepsilon}

\newcounter{marnote}

\begin{document}

\title[Asymptotic analysis for hydrodynamic force]{Asymptotic analysis for hydrodynamic force acting on stiff particles}


\author[Z.W. Zhao]{Zhiwen Zhao}

\address[Z.W. Zhao]{Beijing Computational Science Research Center, Beijing 100193, China.}
\email{zwzhao365@163.com}


\date{\today} 


\maketitle
\begin{abstract}
A three-dimensional mathematical model of a viscous incompressible fluid with two stiff particles is investigated in the near-contact regime. When one of the particles approaches the other motionless particle with prescribed translational and angular velocities, there always appears blow-up of hydrodynamic force exerted on the moving particle. In this paper, we construct explicit singular functions corresponding to the fluid velocity and pressure to establish precise asymptotic formulas for hydrodynamic force with respect to small interparticle distance, which show that its largest singularity is determined by squeeze motion between two particles. Finally, the primal-dual variational principle is employed to give a complete justification for these asymptotics.

\end{abstract}


\section{Introduction}
Many complicated natural phenomena and engineering processes have close relation with suspensions of stiff particles in an incompressible fluid. The investigation on the behavior of suspensions is crucial to its applications in environmental geophysics, chemical engineering, ceramics processing, biotechnology and pharmacology, see e.g. \cite{G1994,S2011,PT2009}. A non-colloidal suspension of neutrally buoyant stiff particles immersed into a Newtonian fluid, where the inertial forces and Brownian motion may be neglected, can be described by Stokes equation. Hydrodynamic force and the effective viscosity are two primary physical quantities of interest in describing the rheological properties of suspensions. In particular, these two quantities will exhibit high singularities in terms of the interparticle distance $\varepsilon$, as the distance $\varepsilon$ goes to zero. Moreover, the contribution to the singularities only comes from thin gaps between neighboring particles.

In this paper, we mainly focus on the blow-up of hydrodynamic force. The problem has a quite long history and has resulted in a long list of literature (for example, \cite{GCB1967,C1974,CB1967,CB1989,J1982,JO1984,G2016}) studying different cases based on shape and number of particles, approximation method, applied boundary conditions, etc. The results provided by these papers show that the largest blow-up rate of hydrodynamic force is of order $O(\varepsilon^{-1})$ in the presence of spherical particles. Particularly in \cite{G2016}, the singular behavior of hydrodynamic force were precisely captured and rigorously justified by combining the polynomial approximation technique and the primal-dual variational principle. Although spherical particle, as the most ideal model of inclusion from the view of theory analysis and numerical computation, is extensively studied in previous literature on the blow-up of hydrodynamic force, it is more practical and essential to study the general non-radial symmetric shape of inclusions frequently in suspensions of nature. Motivated by this fact, Li, Wang and Zhao \cite{LWZ2020} extended the results in \cite{G2016} to the case of two close-to-touching $m$-convex particles with $m\geq2$ (see \eqref{CONVEX001} for the definition) and revealed that hydrodynamic force blows up at the rate of $\varepsilon^{3/m-3}$ and $\varepsilon^{4/m-3}$ in dimensions two and three, respectively. This implies that the singularities of hydrodynamic force increase as the surface convexity of particles weakens. In addition, the shape of particles considered in \cite{C1974} covers the general strictly convex particles, which corresponds to the case of $m=2$ in \cite{LWZ2020}. With regard to the effective viscosity, we refer to \cite{BBP2006,FA1967,G1981,NK1984,BGN2009} and the references therein. In particular, Berlyand, Gorb and Novikov \cite{BGN2009} developed a $fictitious\;fluid\;approach$ to establish the precise asymptotic formulas of the viscous dissipation rate under generic boundary condition in dimension two. Their results also revealed a novel blow-up phenomena of two-dimensional suspensions that the Poiseuille type microflow may result in anomalously strong singularity of the dissipation rate, which is different from the three-dimensional suspensions. In addition, sedimentation in suspensions of stiff particles is another interesting and challenging problem. For this problem a lot of physical phenomena, such as acceleration of sedimentation by the applied shearing in the process of dewatering of waste water sludge in a centrifuge \cite{G2003,RS2001}, are not yet fully understood.

High concentration phenomena caused by the closeness between inclusions or between the inclusion and the external boundary also appears in high-contrast composite materials, which is modeled by the Lam\'{e} system $\mu\Delta\mathbf{u}+(\lambda+\mu)\nabla\nabla\cdot\mathbf{u}=0$. In fact, the Stokes system, which describes the aforementioned concentrated suspensions of rigid particles, is closely related to  the Lam\'{e} system. To be specific, the Lam\'{e} system converges to the Stokes system, as $\lambda\rightarrow\infty$ and $\mu$ is fixed. This fact has been demonstrated by Ammari, Garapon, Kang and Lee \cite{AGKL2008}. The literature on the study of stress concentration occurring in high-contrast composites are very wide and we just mention \cite{KY2019,BLL2015,BLL2017,LZ2019} and the references therein for an interested reader.

The paper is organized as follows. In Section \ref{M001} we set up the problem and state the main results in Theorem \ref{LEADINGRESULTS}, whose proof is given in Section \ref{M002}. Specially, we first carry out a linear decomposition for the original problem \eqref{11.3} according to the elementary kinematic motions of particle, see subsection \ref{EC001}. With regard to these subproblems, we construct the corresponding explicit singular functions for the fluid velocity and pressure to accurately capture the singularities of hydrodynamic force, see subsections \ref{EC002} and \ref{EC003}. It is worth emphasizing that the approximation method presented in subsections \ref{EC002} and \ref{EC003} is greatly different from the polynomial approximation adopted in \cite{G2016}. Finally, the asymptotic results are rigorously justified by utilizing the dual variational principle, see subsection \ref{EC005}.

\section{Formulation of the problem and main results}\label{M001}

Consider a bounded convex domain $D\subset\mathbb{R}^{3}$ which is occupied by an incompressible viscous fluid with viscosity $\mu$. Let $D_{i}\subset D$, $i=1,2$ be two stiff convex particles with $\varepsilon$-apart, where the interparticle distance $\varepsilon$ is a sufficiently small positive constant. Denote $\Omega:=D\setminus\overline{D_{1}\cup D_{2}}$. Assume further that $\partial\Omega$ is of $C^{2,\gamma}$, $0<\gamma<1$ and these two particles are far away from the external boundary $\partial D$. Let the particle $D_{1}$ approach the motionless particle $D_{2}$ with linear velocity $\mathbf{U}$ and angular velocity $\boldsymbol{\omega}$, respectively, written as
\begin{align*}
 \mathbf{U}=
U_{1}\mathbf{e}_{1}+U_{2}\mathbf{e}_{2}+U_{3}\mathbf{e}_{3},\quad \boldsymbol{\omega}=\omega_{1}\mathbf{e}_{1}+\omega_{2}\mathbf{e}_{2}+\omega_{3}\mathbf{e}_{3},
\end{align*}
where $\{\mathbf{e}_{1},\mathbf{e}_{2},\mathbf{e}_{3}\}$ denotes the standard Euclidean basis in $\mathbb{R}^{3}$.

Denote by $\mathbf{u}=(u_{1},u_{2},u_{3})^{T}:D\rightarrow\mathbb{R}^{3}$ and $p:D\rightarrow\mathbb{R}$ the fluid velocity and pressure, respectively. Let $\mathbf{x}_{D_{1}}$ be the center of mass of $D_{1}$. In this paper, we consider the Stokes flow with the following boundary value conditions:
\begin{align}\label{11.3}
\begin{cases}
(\mathrm{a})~\nabla\cdot\sigma[\mathbf{u},p]=\mathbf{0}&~~~~~\mathrm{in}~\Omega,\\
(\mathrm{b})~~~~~~~\nabla\cdot\mathbf{u}=0&~~~~~\mathrm{in}~\Omega,\\
(\mathrm{c})~~~~~~~~~~~~~~\mathbf{u}=\mathbf{U}+\boldsymbol{\omega}\times(\mathbf{x}-\mathbf{x}_{D_{1}})&~~~~~\mathrm{on}~\partial D_{1},\\
(\mathrm{d})~~~~~~~~~~~~~~\mathbf{u}=\mathbf{0}&~~~~~\mathrm{on}~\partial D_{2},\\
(\mathrm{e})~\mathbf{u}=\boldsymbol{\varphi}&~~~~~\mathrm{on}~\partial D,
\end{cases}
\end{align}
where $\sigma[\mathbf{u},p]:=2\mu e(\mathbf{u})-p\mathbb{I}$ represents the Cauchy stress tensor with $e(\mathbf{u}):=\frac{1}{2}(\nabla\mathbf{u}+(\nabla\mathbf{u})^{T})$ and $\mathbb{I}$ denoting the rate of strain tensor and the identity matrix, respectively, $\boldsymbol{\varphi}\in C^{2}(\partial D;\mathbb{R}^{3})$ is the given velocity field verifying the following compatibility condition:
\begin{align}\label{COM001}
\int_{\partial D}\boldsymbol{\varphi}\cdot\mathbf{n}=0.
\end{align}
Here and throughout this paper, $\mathbf{n}$ represents the unit outer normal to the domain. In physics, $(\ref{11.3})(\mathrm{b})$ describes the incompressibility of the Stokes flow. $(\ref{11.3})(\mathrm{c})$ and $(\ref{11.3})(\mathrm{d})$ imply that there appears no-slip between the fluid and the surfaces of particles $D_{1}$ and $D_{2}$. It is worthwhile to point out that the Neumann condition ``$\sigma[\mathbf{u},p]\mathbf{n}=\mathbf{0}$ on $\partial D$" was previously added in \cite{G2016} to let the fluid not flow out of the exterior boundary $\partial D$. However, since the singularity of hydrodynamic force only comes from the narrow channel between two particles, then this condition is actually not essential for the following approximation results to hold. In this paper, we replace it by the Dirichlet condition ``$\mathbf{u}=\boldsymbol{\varphi}$ on $\partial D$" with $\boldsymbol{\varphi}$ satisfying compatibility condition \eqref{COM001}, which is used to ensure the existence and uniqueness of a weak solution to problem \eqref{11.3} (see the detailed proof in \cite{L1959} with a
slight modification). As seen in \cite{T1984}, the Stokes equation is elliptic in the sense of Douglis-Nirenberg. So the general regularity theory established in \cite{ADN1964,S1966} also hold for the Stokes flow.

The principal physical quantities concerned in this paper are {\em hydrodynamic force} exerted on $D_{1}$ and {\em hydrodynamic torque}, respectively, given by
\begin{equation*}
\mathbf{F}=\int_{\partial D_{1}}\sigma[\mathbf{u},p]\mathbf{n}\,dS,
\end{equation*}
and
\begin{equation*}
\mathbf{T}=\int_{\partial D_{1}}(\mathbf{x}-\mathbf{x}_{D_{1}})\times\sigma[\mathbf{u},p]\mathbf{n}\,dS.
\end{equation*}

To state our principal results in a precise manner, we first parameterize the domain. By picking a proper coordinate system, we have
\begin{align*}
\partial(D_{1}+(0',-\varepsilon/2))\cap\partial(D_{2}+(0',\varepsilon/2))=\{0\},
\end{align*}
and
\begin{align*}
D_{i}+(0',(-1)^{i}\varepsilon/2)\subset\{(x',x_{3})\in\mathbb{R}^{3}\,|\,(-1)^{i-1}x_{3}>0\},\quad i=1,2.
\end{align*}
Here and in the following, we denote the two-dimensional variables by adding superscript prime, for example, $0'=(0,0)$ and $x'=(x_{1},x_{2})$.
Let the centers of mass of particles $D_{1}$ and $D_{2}$ be, respectively, located at $\left(0',\pm(\varepsilon/2+R)\right)$ for a small positive constant $R$. Suppose that there exists a small $\varepsilon$-independent constant $0<r<R$ such that the portions of $\partial D_{1}$ and $\partial D_{2}$ around the origin are, respectively, represented by two smooth functions $\pm(\varepsilon/2+h(x'))$ satisfying that
\begin{align}\label{CONVEX001}
h(x')=\kappa|x'|^{m},\quad|x'|\leq r,
\end{align}
where $m\geq2$ and $\kappa$ is a positive constant independent of $\varepsilon$. From the view of the geometry, the curvature of the surfaces of two particles is not equal to zero at points $(0',\pm\varepsilon/2)$ in the case of $m=2$, while it degenerates to be zero for $m>2$. Moreover, the greater the value of the convexity index $m$, the flatter the surfaces of two particles. It is worthwhile to mention that the shape of inclusion considered in condition \eqref{CONVEX001} contains a class of axisymmetric ellipsoids. To be precise, let $\partial D_{1}$ and $\partial D_{2}$ be, respectively, parameterized as
\begin{align*}
|x'|^{m}+|x_{3}\pm(\varepsilon/2+R)|^{m}=R^{m}.
\end{align*}
By Taylor expansion, we have
\begin{align*}
h(x')=\frac{1}{m R^{m-1}}|x'|^{m}+O(|x'|^{2m}),\quad\mathrm{in}\;\Omega_{r}.
\end{align*}
Then this type of axisymmetric ellipsoid corresponds to the case of $\kappa=m^{-1} R^{1-m}$ in condition \eqref{CONVEX001}.

Throughout this paper, for $ij\in\{12,34\}$ and $m\geq2$, denote
\begin{align*}
\Gamma_{ij}^{(m)}&=
\begin{cases}
\frac{1}{m(2\kappa)^{\frac{j}{m}}}\Gamma\left(i-\frac{j}{m}\right)\Gamma\left(\frac{j}{m}\right),&i>\frac{j}{m},\vspace{0.5ex}\\
\frac{1}{m(2\kappa)^{\frac{j}{m}}},&i=\frac{j}{m},
\end{cases}
\end{align*}
where $\Gamma(s)=\int^{+\infty}_{0}x^{s-1}e^{-x}dx$, $s>0$ denotes the Gamma function. In the following, $O(1)$ represents some quantity satisfying that $|O(1)|\leq C$ for some positive constant $C$ independent of $\varepsilon$, which depends only on $\mu$, $R$, $r$, $U_{i}$, $\omega_{i}$, $i=1,2,3$ and $\|\boldsymbol{\varphi}\|_{C^{2}(\partial D)}$.

\begin{theorem}\label{LEADINGRESULTS}
Suppose that $D_{1},D_{2}\subset D\subseteq\mathbb{R}^{3}$ are described as above and condition \eqref{CONVEX001} holds. Let $\mathbf{u}\in H^{1}(D;\mathbb{R}^{3})\cap C^{1}(\overline{\Omega};\mathbb{R}^{3})$ and $p\in L^{2}(D)\cap C^{0}(\overline{\Omega})$ be the solution of problem \eqref{11.3}--\eqref{COM001}. Then for a arbitrarily small $\varepsilon>0$,

$(i)$ if $\boldsymbol{\omega}\neq\mathbf{0}$, then for $m=2$,
\begin{align*}
\mathbf{F}=&-\frac{3\pi\mu}{8\kappa^{2}\varepsilon}\mathbf{U}\cdot\mathbf{e}_{3}-\frac{\pi\mu}{2\kappa}|\ln\varepsilon|\mathbf{U}\cdot(\mathbf{e}_{1}+\mathbf{e}_{2}+\mathbf{e}_{3})-\frac{\pi\mu(10\kappa R-3)}{20\kappa^{2}}|\ln\varepsilon|\boldsymbol{\omega}\times\mathbf{e}_{3}+O(1),\\
\mathbf{T}=&-\frac{\pi\mu R}{2\kappa}|\ln\varepsilon|(\mathbf{U}\times\mathbf{e}_{3}+R\boldsymbol{\omega}\cdot(\mathbf{e}_{1}+\mathbf{e}_{2}))+O(1);
\end{align*}

$(ii)$ if $\boldsymbol{\omega}=\mathbf{0}$, then for $m\geq2$,
\begin{align*}
\mathbf{F}=&-\frac{3\pi\mu\Gamma^{(m)}_{34}}{\varepsilon^{3-4/m}}\mathbf{U}\cdot\mathbf{e}_{3}-\frac{2\pi\mu\Gamma^{(m)}_{12}}{\varepsilon^{1-2/m}}\mathbf{U}\cdot(\mathbf{e}_{1}+\mathbf{e}_{2}+\mathbf{e}_{3})+O(1),\\
\mathbf{T}=&-\frac{2\pi\mu R\Gamma^{(m)}_{12}}{\varepsilon^{1-2/m}}\mathbf{U}\times\mathbf{e}_{3}+O(1).
\end{align*}

\end{theorem}

%
%

\begin{remark}
As shown in Theorem \ref{LEADINGRESULTS}, when $\boldsymbol{\omega}\neq\mathbf{0}$ and $m=2$, the asymptotic expansion of hydrodynamic force shows that its biggest blow-up rate $\varepsilon^{-1}$ is only created by squeeze motion of linear motion, see also \cite{G2016} concerning this conclusion. In the case of $\boldsymbol{\omega}=\mathbf{0}$, we consider the general $m$-convex particles with $m\geq2$ and capture the blow-up rate of order $O(\varepsilon^{4/m-3})$, which implies that the singularity of hydrodynamic force will strengthen as the surface convexity of particles weakens. This fact has also been revealed in previous work \cite{LWZ2020}.

\end{remark}

\section{The proof of Theorem \ref{LEADINGRESULTS}}\label{M002}

\subsection{Linear decomposition}\label{EC001}
Observe that
\begin{equation*}
\mathbf{x}-\mathbf{x}_{D_{1}}=x_{1}\mathbf{e}_{1}+x_{2}\mathbf{e}_{2}+\left(x_{3}-\varepsilon/2-R\right)\mathbf{e}_{3},\quad\mathrm{on}\;\partial D_{1},
\end{equation*}
and thus
\begin{align*}
\mathbf{u}|_{\partial D_{1}}=&\mathbf{U}+\boldsymbol{\omega}\times(\mathbf{x}-\mathbf{x}_{D_{1}})\notag\\
=&\left[U_{1}+\omega_{2}(x_{3}-\varepsilon/2-R)-\omega_{3}x_{2}\right]\mathbf{e}_{1}+\left[U_{2}-\omega_{1}(x_{3}-\varepsilon/2-R)+\omega_{3}x_{1}\right]\mathbf{e}_{2}\notag\\
&+\left(U_{3}+\omega_{1}x_{2}-\omega_{2}x_{1}\right)\mathbf{e}_{3}.
\end{align*}
For $x\in \mathbb{R}^{3}$, let
\begin{align*}
&\boldsymbol{\phi}_{1}=(U_{1}-\omega_{2}R)\mathbf{e}_{1},\;\,\boldsymbol{\phi}_{2}=(U_{2}+\omega_{1}R)\mathbf{e}_{2},\;\,\boldsymbol{\phi}_{3}=U_{3}\mathbf{e}_{3},\;\,\boldsymbol{\phi}_{4}=(\omega_{1}x_{2}-\omega_{2}x_{1})\mathbf{e}_{3},\\
&\boldsymbol{\phi}_{5}=\left(\omega_{2}(x_{3}-\varepsilon/2)-\omega_{3}x_{2}\right)\mathbf{e}_{1}+\left(-\omega_{1}(x_{3}-\varepsilon/2)+\omega_{3}x_{1}\right)\mathbf{e}_{2}.
\end{align*}
Then we have
\begin{align*}
\mathbf{u}=&\sum^{5}_{\alpha=1}\boldsymbol{\phi}_{\alpha},\quad\mathrm{on}\;\partial D_{1}.
\end{align*}
By linearity, the solution $(\mathbf{u},p)$ of original problem \eqref{11.3} can be split as follows:
\begin{align*}
\mathbf{u}=\sum^{5}_{\alpha=1}\mathbf{u}^{(\alpha)},\quad p=\sum^{5}_{\alpha=1}p^{(\alpha)},
\end{align*}
where $(\mathbf{u}^{(\alpha)},p^{(\alpha)})$, $\alpha=1,2,3,4,5$, respectively, verify
\begin{align}\label{Narrow002}
\begin{cases}
\nabla\cdot\sigma[\mathbf{u}^{(\alpha)},p^{(\alpha)}]=\mathbf{0}&\mathrm{in}\;\Omega,\\
\nabla\cdot\mathbf{u}^{(\alpha)}=0&\mathrm{in}\;\Omega,\\
\mathbf{u}^{(\alpha)}=\boldsymbol{\phi}_{\alpha}&\mathrm{on}\;\partial D_{1},\\
\mathbf{u}^{(\alpha)}=\mathbf{0}&\mathrm{on}\;\partial D_{2}\cup\partial D,\\
\end{cases}\quad\alpha=1,2,3,4,
\end{align}
and
\begin{align}\label{Narrow003}
\begin{cases}
\nabla\cdot\sigma[\mathbf{u}^{(5)},p^{(5)}]=\mathbf{0}&\mathrm{in}\;\Omega,\\
\nabla\cdot\mathbf{u}^{(5)}=0&\mathrm{in}\;\Omega,\\
\mathbf{u}^{(5)}=\boldsymbol{\phi}_{5}&\mathrm{on}\;\partial D_{1},\\
\mathbf{u}^{(5)}=\mathbf{0}&\mathrm{on}\;\partial D_{2},\\
\mathbf{u}^{(5)}=\boldsymbol{\varphi}&\mathrm{on}\;\partial D.
\end{cases}
\end{align}
Recall that the classical linear motions of particles comprise of the following three types: parallel translation, shear motion and squeeze motion. According to the aforementioned decomposition, we see that when $\boldsymbol{\omega}=\mathbf{0}$, $\mathbf{u}^{(\alpha)}$, $\alpha=1,2$ correspond to the shear-type motion between two particles, while $\mathbf{u}^{(3)}$ corresponds to the squeeze-type motion between two particles.

For $\alpha=1,2,...,5,$ denote
\begin{align*}
\mathbf{F}^{(\alpha)}=\int_{\partial D_{1}}\sigma[\mathbf{u}^{(\alpha)},p^{(\alpha)}]\mathbf{n}\,dS,\quad\mathbf{T}^{(\alpha)}=\int_{\partial D_{1}}(\mathbf{x}-\mathbf{x}_{D_{1}})\times\sigma[\mathbf{u}^{(\alpha)},p^{(\alpha)}]\mathbf{n}\,dS.
\end{align*}
Then using linearity again, we have
\begin{align}\label{11.17}
\mathbf{F}=\sum\limits^{5}_{\alpha=1}\mathbf{F}^{(\alpha)},\quad\mathbf{T}=\sum\limits^{5}_{\alpha=1}\mathbf{T}^{(\alpha)}.
\end{align}

\subsection{Constructions of the leading terms}\label{EC002}
For $0<t\leq r$, denote the narrow region between two particles by
$$\Omega_{t}=\left\{(x',x_{3})\in\mathbb{R}^{3}\,|~|x_{3}|<\varepsilon/2+\kappa|x'|^{m},~ |x'|<t\right\}.$$
For simplicity, define
\begin{align}\label{DEL010}
\delta:=\delta(x')=\varepsilon+2h(x')=\varepsilon+2\kappa|x'|^{m},\quad\mathrm{in}\;\Omega_{r},
\end{align}
where $\kappa$ is given in \eqref{CONVEX001}. Introduce the following constants:
\begin{align}\label{constants001}
a_{1}=&3,\;\, a_{2}=-2,\;\,b_{1}=-\frac{12}{5},\;\,b_{2}=\frac{3}{10\kappa},\;\,b_{3}=\frac{16}{5},\;\,b_{4}=\frac{3}{5\kappa}.
\end{align}

Since the contribution to the singularities of hydrodynamic force only comes from the thin gap between two particles, then it is the key to give the explicit singular functions corresponding to the solution $(\mathbf{u}^{(\alpha)},p^{(\alpha)})$ in $\Omega_{r}$ for the purpose of accurately calculating $\mathbf{F}$ and $\mathbf{T}$ in \eqref{11.17}, $\alpha=1,2,...,5.$ Introduce a family of auxiliary functions $\bar{\mathbf{u}}^{(\alpha)}\in C^{2,\gamma}(\Omega;\mathbb{R}^{3})$, $\alpha=1,2,...,5$, such that $\|\bar{\mathbf{u}}^{(\alpha)}\|_{C^{2,\gamma}(\Omega\setminus\Omega_{r})}\leq C$, $\nabla\cdot\bar{\mathbf{u}}^{(\alpha)}=0$ in $\Omega_{r}$,
\begin{align}\label{WAE001}
\begin{cases}
\bar{\mathbf{u}}^{(\alpha)}=\boldsymbol{\phi}_{\alpha}&\mathrm{on}\;\partial D_{1},\\
\bar{\mathbf{u}}^{(\alpha)}=\mathbf{0}&\mathrm{on}\;\partial D_{2}\cup\partial D,\\
\end{cases}\;\,\alpha=1,2,3,4,\quad
\begin{cases}
\mathbf{u}^{(5)}=\boldsymbol{\phi}_{5}&\mathrm{on}\;\partial D_{1},\\
\mathbf{u}^{(5)}=\mathbf{0}&\mathrm{on}\;\partial D_{2},\\
\mathbf{u}^{(5)}=\boldsymbol{\varphi}&\mathrm{on}\;\partial D,
\end{cases}
\end{align}
and
\begin{align}\label{zzwz002}
\bar{\mathbf{u}}^{(\alpha)}=&\boldsymbol{\phi}_{\alpha}\left(\frac{1}{2}+\mathfrak{G}\right)+\left(\mathfrak{G}^{2}-\frac{1}{4}\right)\boldsymbol{\mathcal{F}}_{\alpha},\quad \mathrm{in}\;\Omega_{r},
\end{align}
where
\begin{align}\label{deta}
\mathfrak{G}:=\mathfrak{G}(x)=\frac{x_{3}}{\delta},
\end{align}
and $\boldsymbol{\mathcal{F}}_{\alpha}=(\mathcal{F}^{(\alpha)}_{1},\mathcal{F}^{(\alpha)}_{2},\mathcal{F}^{(\alpha)}_{3})$, $\alpha=1,2,...,5$ are, respectively, given by
\begin{align}\label{QLA001}
\boldsymbol{\mathcal{F}}_{\alpha}=&
\begin{cases}
\frac{(U_{1}-\omega_{2}R)\partial_{1}\delta}{2}\mathbf{e}_{3},&\alpha=1,\,m\geq2,\\
\frac{(U_{2}+\omega_{1}R)\partial_{2}\delta}{2}\mathbf{e}_{3},&\alpha=2,\,m\geq2,\\
U_{3}\Big[\sum\limits^{2}_{i=1}\frac{a_{1}x_{i}}{\delta}\mathbf{e}_{i}+\frac{x_{3}}{\delta}\big(\frac{a_{1}(x'\cdot\nabla_{x'}\delta)}{\delta}+a_{2}\big)\mathbf{e}_{3}\Big],&\alpha=3,\,m\geq2,
\end{cases}
\end{align}
and, for $\alpha=4$ and $m=2$,
\begin{align}\label{DANZW001}
\boldsymbol{\mathcal{F}}_{4}=&\sum\limits^{2}_{i=1}\frac{b_{1}x_{i}(-\omega_{1}x_{2}+\omega_{2}x_{1})}{\delta}\mathbf{e}_{i}+b_{2}(-\omega_{1}\mathbf{e}_{2}+\omega_{2}\mathbf{e}_{1})\notag\\
&+(-\omega_{1}x_{2}+\omega_{2}x_{1})\Big(\frac{b_{1}(x'\cdot\nabla_{x'}\delta)}{\delta}+b_{3}\Big)\mathfrak{G}\mathbf{e}_{3},
\end{align}
and, for $\alpha=5$ and $m=2$,
\begin{align*}
\boldsymbol{\mathcal{F}}_{5}=&\Big[\frac{1}{2}(\omega_{2}(x_{3}-\varepsilon/2)-\omega_{3}x_{2})\partial_{1}\delta+\frac{1}{2}(-\omega_{1}(x_{3}-\varepsilon/2)+\omega_{3}x_{1})\partial_{2}\delta\notag\\
&\,+\frac{1}{4}\delta(\omega_{1}\partial_{2}\delta-\omega_{2}\partial_{1}\delta)\Big]\mathbf{e}_{3}.
\end{align*}

Remark that for $\alpha=1,2,3,4$, the first part $\boldsymbol{\phi}_{\alpha}\left(\frac{1}{2}+\mathfrak{G}\right)$ constructed in \eqref{zzwz002} is linear in $x_{3}$ in the neck $\Omega_{r}$ with the values $\boldsymbol{\phi}_{\alpha}$ and $\mathbf{0}$ on the top and bottom boundaries $\Sigma^{\pm}_{r}:=\left\{(x',x_{3})\in\mathbb{R}^{3}\,|~ x_{3}=\pm(\varepsilon/2+h(x')),~ |x'|<r\right\}$. This type of test function is called the Keller-type function, which was first introduced in \cite{K1963} to derive the effective conductivity. Since then, it has been extensively used to study the stress concentration for the elasticity problem arising from high-contrast composites (see, for example, \cite{BK2001,BLL2015,BLL2017} and references therein). The second part $\left(\mathfrak{G}^{2}-\frac{1}{4}\right)\boldsymbol{\mathcal{F}}_{\alpha}$ is called the correction term, which belongs to the lower order singular term and can be determined by using the incompressibility of Stokes flow.

Although the Keller-type test function is a key component in the constructions of the leading terms both for the elasticity problem and the Stokes problem, there is a significant difference that Stokes equation contains the pressure. So in order to apply the dual variational principle to complete the justification for the following approximation results, we need to construct the corresponding test functions for the pressure to eliminate the biggest singular terms contained in $\mu\partial_{nn}\bar{\mathbf{u}}^{(\alpha)}$, $\alpha=1,2,...,5$. For this purpose, we introduce scalar test functions $\bar{p}^{(\alpha)}\in C^{1,\gamma}(\Omega)$, $\alpha=1,2,...,5$ satisfying that $\|\bar{p}^{(\alpha)}\|_{C^{1,\gamma}(\Omega\setminus\Omega_{r})}\leq C$, and for $x\in\Omega_{r}$,
\begin{align}\label{PMAIN002}
\bar{p}^{(\alpha)}=&
\begin{cases}
\frac{\mu(U_{1}-\omega_{2}R) x_{3}\partial_{1}\delta}{\delta^{2}},&\alpha=1,\,m\geq2,\vspace{0.5ex}\\
\frac{\mu(U_{2}+\omega_{1}R) x_{3}\partial_{2}\delta}{\delta^{2}},&\alpha=2,\,m\geq2,\vspace{0.5ex}\\
\mu U_{3}\Big[\frac{3 x_{3}^{2}}{\delta^{3}}\big(\frac{a_{1}(x'\cdot\nabla_{x'}\delta)}{\delta}+a_{2}\big)+ a_{1}\int^{|x'|^{2}}_{r^{2}}\frac{1}{(\varepsilon+2\kappa t^{\frac{m}{2}})^{3}}dt\Big],&\alpha=3,\,m\geq2,\vspace{0.5ex}\\
\frac{\mu (\omega_{2}x_{1}-\omega_{1}x_{2})}{\delta^{2}}\Big[\frac{3x_{3}^{2}}{\delta}\big(\frac{b_{1}(x'\cdot\nabla_{x'}\delta)}{\delta}+b_{3}\big)+b_{4}\Big],&\alpha=4,\,m=2,\vspace{0.5ex}\\
0,&\alpha=5,\,m=2,
\end{cases}
\end{align}
where $a_{i}$, $i=1,2,$ and $b_{j},$ $j=1,3,4$ are defined by \eqref{constants001}. Combining these singular functions, we obtain

$(1)$ for $i=1,2$,
\begin{align}\label{DNZ001}
&\mu\partial_{33}\bar{\mathbf{u}}^{(\alpha)}_{i}-\partial_{i}\bar{p}^{(\alpha)}\notag\\
&=
\begin{cases}
-\mu(U_{1}-\omega_{2}R)x_{3}\partial_{i}(\delta^{-2}\partial_{1}\delta),&\alpha=1,\,m\geq2,\vspace{0.5ex}\\
-\mu(U_{2}+\omega_{1}R)x_{3}\partial_{i}(\delta^{-2}\partial_{2}\delta),&\alpha=2,\,m\geq2,\vspace{0.5ex}\\
-3\mu U_{3}x_{3}^{2}\partial_{i}\big[\delta^{-3}\big(\frac{a_{1}(x'\cdot\nabla_{x'}\delta)}{\delta}+a_{2}\big)\big],&\alpha=3,\,m\geq2,\vspace{0.5ex}\\
-3\mu x_{3}^{2}\partial_{i}\big[(\omega_{2}x_{1}-\omega_{1}x_{2})\delta^{-3}\big(\frac{b_{1}(x'\cdot\nabla_{x'}\delta)}{\delta}+b_{3}\big)\big],&\alpha=4,\,m=2,\vspace{0.5ex}\\
0,&\alpha=5,\,m=2;
\end{cases}
\end{align}

$(2)$ for $i=3$,
\begin{align}\label{DNZ002}
&\mu\partial_{33}\bar{\mathbf{u}}_{3}^{(\alpha)}-\partial_{3}\bar{p}^{(\alpha)}=0,\quad\alpha=1,2,3,4,
\end{align}
and
\begin{align}\label{DNZ003}
\mu\partial_{33}\bar{\mathbf{u}}_{3}^{(5)}-\partial_{3}\bar{p}^{(5)}=&\mu\delta^{-2}[(\omega_{2}(x_{3}-\varepsilon/2)-\omega_{3}x_{2})\partial_{1}\delta+(-\omega_{1}(x_{3}-\varepsilon/2)+\omega_{3}x_{1})\partial_{2}\delta]\notag\\
&+\frac{\mu}{2\delta}(\omega_{1}\partial_{2}\delta-\omega_{2}\partial_{1}\delta).
\end{align}

We here would like to point out that the aforementioned singular functions correspond to the following two cases: when $\boldsymbol{\omega}\neq\mathbf{0}$, we consider $m=2$; when $\boldsymbol{\omega}=\mathbf{0}$, we consider $m\geq2$. Inspired by similar constructions of the auxiliary functions in \cite{LX2022}, we can use the method of undetermined coefficients to solve these singular functions. We now present the solution procedure for finding the correction terms in \eqref{QLA001}. Take $\bar{\mathbf{u}}^{(\alpha)}$, $\alpha=3,4$ for instance. Other cases are the same.

$(i)$ Consider the case when $\alpha=3$ and $m\geq2$. Utilizing \eqref{DEL010}, \eqref{deta} and the incompressible condition, it follows from a direct computation that
\begin{align*}
0=\nabla\cdot\bar{\mathbf{u}}^{(3)}=U_{3}\left[\left(1-\frac{2a_{1}+a_{2}}{4}\right)\partial_{3}\mathfrak{G}+(2a_{1}+3a_{2})\mathfrak{G}^{2}\partial_{3}\mathfrak{G}\right],
\end{align*}
which implies that $a_{1}=3$ and $a_{2}=-2$.

$(ii)$ Consider the case when $\alpha=4$ and $m=2$. Note that
\begin{align*}
\partial_{i}\delta=4\kappa x_{i},\;\, \partial_{i}\mathfrak{G}=-4\kappa x_{i}\mathfrak{G}\partial_{3}\mathfrak{G},\quad i=1,2,
\end{align*}
which, together with \eqref{zzwz002} and \eqref{DANZW001}, reads that
\begin{align*}
\nabla\cdot\bar{\mathbf{u}}^{(4)}=&(\omega_{1}x_{2}-\omega_{2}x_{1})\left[-\left(1+\frac{3b_{1}+b_{3}}{4}\right)\partial_{3}\mathfrak{G}+(3b_{1}+3b_{3}-8\kappa b_{2})\mathfrak{G}^{2}\partial_{3}\mathfrak{G}\right].
\end{align*}
Applying the incompressible condition to $\bar{\mathbf{u}}^{(4)}$, we obtain
\begin{align*}
\begin{cases}
1+\frac{3b_{1}+b_{3}}{4}=0,\\
3b_{1}+3b_{3}-8\kappa b_{2}=0,
\end{cases}
\end{align*}
which yields that
\begin{align}\label{AGZW001}
b_{1}=-\frac{6+4\kappa b_{2}}{n},\;\, b_{3}=2+4\kappa b_{2}.
\end{align}
For later justification, we now use $\partial_{i}\bar{p}^{(4)}$, $i=1,2$ to offset the greatest singular terms of order $O(\delta^{-2})$ in $\mu\partial_{33}\bar{\mathbf{u}}^{(4)}_{i}$, $i=1,2.$ That is, let
\begin{align*}
\begin{cases}
\frac{2\mu}{\delta^{2}}\Big(\frac{b_{1}(\omega_{2}x_{1}^{2}-\omega_{1}x_{2}^{2})}{\delta}+b_{2}\Big)=\mu b_{4}[\omega_{2}\partial_{1}(x_{1}\delta^{-2})-\omega_{1}\partial_{2}(x_{2}\delta^{-2})],\vspace{0.5ex}\\
\frac{2\mu b_{1}(\omega_{2}-\omega_{1})x_{1}x_{2}}{\delta^{3}}=\mu b_{4}(\omega_{2}x_{1}\partial_{2}\delta^{-2}-\omega_{1}x_{2}\partial_{1}\delta^{-2}).
\end{cases}
\end{align*}
This leads to that
\begin{align*}
\begin{cases}
(2b_{2}-b_{4})(\omega_{2}-\omega_{1})+(2b_{1}+8\kappa b_{4})(\omega_{2}x_{1}^{2}-\omega_{1}x_{2}^{2})\delta^{-1}=0,\\
(2b_{1}+8\kappa b_{4})(\omega_{2}-\omega_{1})x_{1}x_{2}=0.
\end{cases}
\end{align*}
Hence,
\begin{align*}
b_{1}=-4\kappa b_{4},\;\,b_{2}=\frac{b_{4}}{2},
\end{align*}
which, in combination with \eqref{AGZW001}, shows that
\begin{align*}
b_{1}=-\frac{12}{5},\;\,b_{2}=\frac{3}{10\kappa},\;\,b_{3}=\frac{16}{5},\;\,b_{4}=\frac{3}{5\kappa}.
\end{align*}

For readers' convenience, we sum up the method to construct these above-mentioned singular functions. To begin with, we pick the Keller-type test function as the dominant term of the velocity and then use the incompressibility of Stokes flow to find the corresponding correction term. The singular function for the pressure is subsequently chosen to ensure that the pair of singular functions corresponding to the velocity and pressure approximately satisfies the Stokes equation in the sense of no large singular terms appearing in the remainder.

For simplicity, for $ij\in\{12,34\}$ and $m\geq2$, write
\begin{align*}
\rho^{(m)}_{ij}(\varepsilon)=&
\begin{cases}
\frac{1}{\varepsilon^{i-\frac{j}{m}}},&i>\frac{j}{m},\vspace{0.5ex}\\
|\ln\varepsilon|,&i=\frac{j}{m}.
\end{cases}
\end{align*}

Based on the explicit main terms constructed above, we are now ready to give precise calculations for the force and torque.
\subsection{Asymptotics of the force and torque.}\label{EC003}
For $\alpha=1,2,...,5,$ denote by $\overline{\mathbf{F}}^{(\alpha)}$ and $\overline{\mathbf{T}}^{(\alpha)}$ the corresponding approximations for $\mathbf{F}^{(\alpha)}$ and $\mathbf{T}^{(\alpha)}$, respectively. Observe that every constituent part of the main terms in \eqref{zzwz002} and \eqref{PMAIN002} exhibits explicit parity with respect to $x_{i}$ in $\Omega_{r}$, $i=1,2.$ Utilizing this fact, it follows from a direct calculation that

$(i)$ for $\alpha=1$ and $m\geq2$,
\begin{align}\label{ASYM001}
\begin{cases}
\overline{\mathbf{F}}^{(1)}=-2\pi\mu(U_{1}-\omega_{2}R)\Gamma_{12}^{(m)}\rho^{(m)}_{12}(\varepsilon)\,\mathbf{e}_{1}+O(1),\\
\overline{\mathbf{T}}^{(1)}=2\pi\mu R(U_{1}-\omega_{2}R)\Gamma_{12}^{(m)}\rho^{(m)}_{12}(\varepsilon)\,\mathbf{e}_{2}+O(1);
\end{cases}
\end{align}

$(ii)$ for $\alpha=2$ and $m\geq2$,
\begin{align}\label{ASYM002}
\begin{cases}
\overline{\mathbf{F}}^{(2)}=-2\pi\mu(U_{2}+\omega_{1}R)\Gamma_{12}^{(m)}\rho^{(m)}_{12}(\varepsilon)\,\mathbf{e}_{2}+O(1),\\
\overline{\mathbf{T}}^{(2)}=-2\pi\mu R(U_{2}+\omega_{1}R)\Gamma_{12}^{(m)}\rho^{(m)}_{12}(\varepsilon)\,\mathbf{e}_{1}+O(1);
\end{cases}
\end{align}

$(iii)$ for $\alpha=3$ and $m\geq2$,
\begin{align}\label{ASYM003}
\overline{\mathbf{F}}^{(3)}&=-\pi\mu U_{3}\left(2\Gamma^{(m)}_{12}\rho^{(m)}_{12}(\varepsilon)+3\Gamma^{(m)}_{34}\rho^{(m)}_{34}\right)\,\mathbf{e}_{3}+O(1),\quad\overline{\mathbf{T}}^{(3)}=\mathbf{0}+O(1);
\end{align}

$(iv)$ for $\alpha=4$ and $m=2$,
\begin{align}\label{ASYM005}
\begin{cases}
\overline{\mathbf{F}}^{(4)}=-\frac{3\pi\mu\omega_{2}}{5\kappa}\Gamma^{(2)}_{12}\rho^{(2)}_{12}(\varepsilon)\,\mathbf{e}_{1}+\frac{3\pi\mu\omega_{1}}{5\kappa}\Gamma^{(2)}_{12}\rho^{(2)}_{12}(\varepsilon)\,\mathbf{e}_{2}+O(1),\\
\overline{\mathbf{T}}^{(4)}=\mathbf{0}+O(1);
\end{cases}
\end{align}

$(v)$ for $\alpha=5$ and $m=2$,
\begin{align}\label{ASYM006}
\overline{\mathbf{F}}^{(5)}=\mathbf{0}+O(1),\quad \overline{\mathbf{T}}^{(5)}=\mathbf{0}+O(1).
\end{align}

In the next subsection, these approximation results obtained in \eqref{ASYM001}--\eqref{ASYM006} will be justified by taking advantage of the dual variational principle, which was previously presented in \cite{G2016}. By contrast with the justification in \cite{G2016}, there exist some differences, since we consider the Dirichlet-type condition on the external boundary for problem \eqref{11.3} but not the Neumann-type condition.
\subsection{Justification}\label{EC005}
To begin with, by the variational argument, we know that the solution $\mathbf{u}$ of \eqref{11.3} minimizes the energy functional as follows:
\begin{align}\label{5.3}
\mathbf{u}&=\mathrm{arg}\min_{\mathbf{v}\in\mathcal{A}_{\Omega}}I_{\Omega}[\mathbf{v}],~~I_{\Omega}[\mathbf{v}]=\mu\int_{\Omega}(e(\mathbf{v}),e(\mathbf{v}))\,dx,
\end{align}
where $\mathcal{A}_{\Omega}=\{\mathbf{v}\in H^{1}(\Omega)|~\mathbf{v}\text{ satisfies conditions (b)--(e) in \eqref{11.3}}\}.$ Here and in the following, for any two $3\times3$ matrices $\mathbb{A}=(a_{ij})$ and $\mathbb{B}=(b_{ij})$, $$(\mathbb{A},\mathbb{B}):=\mathrm{tr}(\mathbb{A}\mathbb{B})=\sum\limits^{3}_{i,j=1}a_{ij}b_{ij}.$$
Denote
\begin{align*}
\overline{\mathbf{F}}:=\sum^{5}_{\alpha=1}\overline{\mathbf{F}}^{(\alpha)},\quad\overline{\mathbf{T}}:=\sum^{5}_{\alpha=1}\overline{\mathbf{T}}^{(\alpha)},
\end{align*}
where $\overline{\mathbf{F}}^{(\alpha)}$ and $\overline{\mathbf{F}}^{(\alpha)}$, $\alpha=1,2,...,5$ are given by \eqref{ASYM001}--\eqref{ASYM006}. In light of \eqref{11.17}, we see that $\overline{\mathbf{F}}$ and $\overline{\mathbf{T}}$ are actually the approximations corresponding to $\mathbf{F}$ and $\mathbf{T}$, respectively. Define
\begin{align}\label{MAINZ001}
\bar{\mathbf{u}}=\sum^{5}_{\alpha=1}\bar{\mathbf{u}}^{(\alpha)},
\end{align}
where $\bar{\mathbf{u}}^{(\alpha)}$, $\alpha=1,2,...,5$ are defined by \eqref{WAE001}--\eqref{zzwz002}. Then applying integration by parts for $I_{\Omega}[\mathbf{u}]$ and $I_{\Omega}[\bar{\mathbf{u}}]$ with $\mathbf{u}$ being the solution of $\eqref{11.3}$ and $\bar{\mathbf{u}}$ given by \eqref{MAINZ001}, it follows from \eqref{WAE001}--\eqref{zzwz002} and \eqref{DNZ001}--\eqref{DNZ003} that
\begin{align}\label{ANTQ001}
\mathbf{U}\cdot(\overline{\mathbf{F}}-\mathbf{F})+\boldsymbol{\omega}\cdot(\overline{\mathbf{T}}-\mathbf{T})
=2(I_{\Omega}[\mathbf{u}]-I_{\Omega}[\bar{\mathbf{u}}])+O(1).
\end{align}
By the dual variational principle corresponding to $(\ref{5.3})$, we know that the Cauchy stress tensor $\sigma[\mathbf{u},p]$ maximizes the following functional:
\begin{align}\label{11.78}
\sigma[\mathbf{u},p]&=\mathrm{arg}\max_{\mathbb{S}\in\mathcal{A}^{\ast}_{\Omega}}I^{\ast}_{\Omega}[\mathbb{S}],\;I^{\ast}_{\Omega}[\mathbb{S}]=\int_{\partial\Omega}\mathbf{u}\cdot\mathbb{S}\mathbf{n}\,dS-\frac{1}{4\mu}\int_{\Omega}\left(\mathrm{tr}\mathbb{S}^{2}-\frac{(\mathrm{tr}\mathbb{S})^{2}}{3}\right),
\end{align}
where $(\mathbf{u},p)$ is the solution of \eqref{11.3},
\begin{align*}
\mathcal{A}^{\ast}_{\Omega}&=\big\{\mathbb{S}\in\mathbb{R}^{3\times3}|~\mathbb{S}=\mathbb{S}^{T},~\nabla\cdot\mathbb{S}=\mathbf{0}~\mathrm{in}~\Omega,~S_{ij}\in L^{2}(\Omega),\,i,j=1,2,3\big\}.
\end{align*}
According to \eqref{5.3} and \eqref{11.78}, we obtain that for any $\mathbb{S}\in\mathcal{A}^{\ast}_{\Omega}$,
\begin{align}\label{CRIT001}
|I_{\Omega}[\mathbf{u}]-I^{\ast}_{\Omega}[\mathbb{S}]|\leq |I_{\Omega}[\bar{\mathbf{u}}]-I^{\ast}_{\Omega}[\mathbb{S}]|.
\end{align}
Combining \eqref{ANTQ001} and \eqref{CRIT001}, the proof is reduced to finding a test tensor $\overline{\mathbb{S}}\in\mathcal{A}^{\ast}_{\Omega}$ such that
\begin{align*}
\mathrm{Err}:=&I_{\Omega}[\bar{\mathbf{u}}]-I^{\ast}_{\Omega}[\overline{\mathbb{S}}]\notag\\
=&\mu\int_{\Omega}(e(\bar{\mathbf{u}}),e(\bar{\mathbf{u}}))dx-\int_{\partial\Omega}\bar{\mathbf{u}}\cdot\overline{\mathbb{S}}\mathbf{n}~dS+\frac{1}{4\mu}\int_{\Omega}\left(\mathrm{tr}\overline{\mathbb{S}}^{2}-\frac{(\mathrm{tr}\overline{\mathbb{S}})^{2}}{3}\right)dx\notag\\
=&O(1),
\end{align*}
where we used the fact that $\bar{\mathbf{u}}=\mathbf{u}$ on $\partial\Omega$.

Pick
\begin{align*}
\overline{\mathbb{S}}=
\begin{cases}
\sum\limits^{5}_{\alpha=1}\overline{\mathbb{S}}^{(\alpha)},&~~~~~\Omega_{r},\\
\mathbb{O},&~~~~~\Omega\setminus\Omega_{r},
\end{cases}
\end{align*}
where $\mathbb{O}$ denotes null matrix, $\overline{\mathbb{S}}^{(\alpha)}\in\mathbb{R}^{3\times3}$, $\alpha=1,2,...,5$, correspond to the dual variational formulation to problem \eqref{5.3} with $\Omega$ replaced by $\Omega_{r}$, given by
\begin{align*}
\sigma[\mathbf{u}^{(\alpha)},p^{(\alpha)}]&=\mathrm{arg}\max_{\mathbb{S}\in\mathcal{A}^{\ast}_{\Omega_{r}}}I^{\ast}_{\Omega_{r}}[\mathbb{S}],\;I^{\ast}_{\Omega_{r}}[\mathbb{S}]=\int_{\partial\Omega_{r}}\mathbf{u}\cdot\mathbb{S}\mathbf{n}\,dS-\frac{1}{4\mu}\int_{\Omega_{r}}\left(\mathrm{tr}\mathbb{S}^{2}-\frac{(\mathrm{tr}\mathbb{S})^{2}}{3}\right),
\end{align*}
where $$\mathcal{A}^{\ast}_{\Omega_{r}}=\big\{\mathbb{S}\in\mathbb{R}^{3\times3}|~\mathbb{S}=\mathbb{S}^{T},~\nabla\cdot\mathbb{S}=\mathbf{0}~\mathrm{in}~\Omega_{r},~S_{ij}\in L^{2}(\Omega_{r}),\,i,j=1,2,3\big\}.$$
In light of the definitions of $\bar{\mathbf{u}}$ and $\overline{\mathbb{S}}$, it follows from integration by parts and the standard interior and boundary estimates for Stokes equation that
\begin{align}\label{11.85}
\mathrm{Err}=&I_{\Omega_{r}}[\bar{\mathbf{u}}]-I^{\ast}_{\Omega_{r}}[\overline{\mathbb{S}}]+O(1)\notag\\
=&\mu\int_{\Omega_{r}}\mathrm{tr}\left[e(\bar{\mathbf{u}})-\frac{1}{2\mu}\left(\overline{\mathbb{S}}-\frac{\mathrm{tr}\overline{\mathbb{S}}}{3}\mathbb{I}\right)\right]^{2}dx+O(1)\notag\\
=&:\sum\limits^{5}_{\alpha=1}\ell[\alpha,\alpha]+2\sum_{\alpha<\beta}^{5}\ell[\alpha,\beta]+O(1),
\end{align}
where $\mathbb{I}\in\mathbb{R}^{3\times3}$ denotes the unit tensor,
$$\ell[\alpha,\beta]=\mu\int_{\Omega_{r}}\left(e(\bar{\mathbf{u}}^{(\alpha)})-\frac{1}{2\mu}\Big(\overline{\mathbb{S}}^{(\alpha)}-\frac{\mathrm{tr}\overline{\mathbb{S}}^{(\alpha)}}{3}\mathbb{I}\Big),e(\bar{\mathbf{u}}^{(\beta)})-\frac{1}{2\mu}\Big(\overline{\mathbb{S}}^{(\beta)}-\frac{\mathrm{tr}\overline{\mathbb{S}}^{(\beta)}}{3}\mathbb{I}\Big)\right).$$

Then the next goal is to pick appropriate test tensors $\overline{\mathbb{S}}^{(\alpha)}\in\mathcal{A}^{\ast}_{\Omega_{r}}$ such that $\ell[\alpha,\beta]=O(1)$, $\alpha,\beta=1,2,...,5$. Due to the fact that $\overline{\mathbf{u}}^{(5)}$ contributes to no singularity in the asymptotics of $\overline{\mathbf{F}}$ and $\overline{\mathbf{T}}$, we take
\begin{equation}\label{11.90}
\mathbb{S}^{(5)}=\mathbb{O},\quad\mathrm{in}\;\Omega_{r},
\end{equation}
which yields that
\begin{align}\label{11.91}
\ell[5,5]=O(1).
\end{align}
For $\alpha=1,2$, pick
\begin{align}\label{11.92}
\overline{\mathbb{S}}^{(1)}&=\mu(U_{1}-\omega_{2}R)
\left(
\begin{array}{ccc}
          0 & 0 & \delta^{-1} \\
          0 & 0 &  0     \\
         \delta^{-1} & 0 & -x_{3}\partial_{1}\delta^{-1}
\end{array}
 \right),\quad\mathrm{in}~\Omega_{r},
\end{align}
and
\begin{align}\label{11.95}
\mathbb{S}^{(2)}&=\mu(U_{2}+\omega_{1}R)
\left(
\begin{array}{ccc}
          0 & 0 &  0      \\
          0 & 0 & \delta^{-1} \\
          0 & \delta^{-1} &-x_{3}\partial_{2}\delta^{-1}
\end{array}
 \right),\quad\mathrm{in}\;\Omega_{r},
\end{align}
where $\delta$ is defined by \eqref{DEL010}. Then it follows from a straightforward computation that
\begin{align}\label{11.93}
\ell[\alpha,\alpha]=O(1),\quad\alpha=1,2.
\end{align}

For $\alpha=3,4$, choose $\overline{\mathbb{S}}^{(\alpha)}=(\overline{\mathbb{S}}_{kl}^{(\alpha)})_{3\times3}$ such that
\begin{align}\label{11.97}
\overline{\mathbb{S}}_{kl}^{(\alpha)}=
\begin{cases}
2\mu\partial_{k}\overline{\mathbf{u}}_{k}^{(\alpha)}-\overline{p}^{(\alpha)}-\int^{x_{k}}_{0}(\mu\Delta\bar{\mathbf{u}}^{(\alpha)}_{k}-\partial_{k}\bar{p}^{(\alpha)})dx_{k},&k=l,\\
\mu(\partial_{l}\overline{\mathbf{u}}^{(\alpha)}_{k}+\partial_{k}\overline{\mathbf{u}}_{l}^{(\alpha)}),&k\neq l,
\end{cases}
\end{align}
where $\bar{\mathbf{u}}^{(\alpha)}=(\bar{\mathbf{u}}^{(\alpha)}_{1},\bar{\mathbf{u}}^{(\alpha)}_{2},\bar{\mathbf{u}}^{(\alpha)}_{3})$ and $\bar{p}^{(\alpha)}$, $\alpha=3,4$ are, respectively, given by \eqref{zzwz002} and \eqref{PMAIN002}. Then we infer from \eqref{DNZ001}--\eqref{DNZ003} that
\begin{align}\label{11.99}
\ell[\alpha,\alpha]=O(1),\quad \alpha=3,4.
\end{align}

Finally, using the symmetry of integral domain and the parity of integrand, we deduce from \eqref{11.90}, \eqref{11.92}, \eqref{11.95} and \eqref{11.97} that for $\alpha,\beta=1,2,...,5$,
\begin{align}\label{11.100}
\ell[\alpha,\beta]&=O(1),\quad \alpha<\beta.
\end{align}
Substituting \eqref{11.91}, \eqref{11.93} and \eqref{11.99}--\eqref{11.100} into \eqref{11.85}, we derive
\begin{align*}
\mathrm{Err}=O(1).
\end{align*}
Consequently, the proof of Theorem \ref{LEADINGRESULTS} is complete.

\noindent{\bf{\large Acknowledgements.}} The author would like to thank Prof. C.X. Miao for his constant encouragement and useful discussions. The author was partially supported by CPSF (2021M700358).

\end{document}